\documentclass[a4paper,draft]{article}  %,12pt
\usepackage{fullpage}
\usepackage{amsthm}
\usepackage{enumerate}
\usepackage{amsxtra} %for \text
\usepackage{amssymb} %for \mathbb

\newtheorem{Thm}{Theorem}[section] %[subsection]when subsection changes then the counter restarts...
\newtheorem{Def}[Thm]{Definition}
\newtheorem{Lem}[Thm]{Lemma}
\newtheorem{Pro}[Thm]{Proposition}
\newtheorem{Cor}[Thm]{Corollary}

\theoremstyle{definition}
\newtheorem{Rem}[Thm]{Remark}

\numberwithin{equation}{section}

\newcommand{\One}{\textbf{1}}

\newcommand{\E}{\ensuremath{\mathbb{E}}}
\newcommand{\R}{\ensuremath{\mathbb{R}}}

\newcommand{\Prob}{\ensuremath{\mathbb{P}}}
\newcommand{\Q}{\ensuremath{\mathbb{Q}}}

\newcommand{\N}{\ensuremath{\mathbb{N}}}

\newcommand{\Fskript}{\ensuremath{\mathcal{F}}}

\newcommand{\Bskript}{\ensuremath{\mathcal{B}}}

\newcommand{\Borel}{\ensuremath{\mathcal{B}}}

\newcommand{\address}[1]{\begin{center}{\small \textit{#1}}\end{center}}
\newcommand{\sep}{, }
\newenvironment{frontmatter}{}{}
\newenvironment{keyword}{\textit{Keywords: }{}}

\begin{document}
\begin{frontmatter}

\title{An Overshoot Approach to \\Recurrence and Transience of Markov Processes}

\author{Bj\"orn B\"ottcher}
\date{July, 2010}
\maketitle
\address{Fakult{\"a}t Mathematik und Naturwissenschaften, 
Institut f{\"u}r mathematische Stochastik,\\
01062 Dresden,
Germany,
{bjoern.boettcher at tu-dresden.de}}

\begin{abstract}
% The paper is motivated by the open problem of stability of non symmetric Feller processes.

We develop criteria for recurrence and transience of one-dimensional Markov processes which have jumps and oscillate between $+\infty$ and $-\infty$. The conditions are based on a Markov chain which only consists of jumps (overshoots) of the process into complementary parts of the state space.

In particular we show that a stable-like process with generator $-(-\Delta)^{\alpha(x)/2}$ such that $\alpha(x)=\alpha$ for $x<-R$ and $\alpha(x)=\beta$ for $x>R$ for some $R>0$ and $\alpha,\beta\in(0,2)$ is transient if and only if $\alpha+\beta<2$, otherwise it is recurrent.

As a special case this yields a new proof for the recurrence, point recurrence and transience of symmetric $\alpha$-stable processes. 

\end{abstract}


\end{frontmatter}

\section{Introduction}

The recurrence and transience of Markov process has been studied by various authors and various techniques, there is the potential theoretic approach (see Getoor \cite{Geto1980} for a unification of the criteria) and the Markov chain approach by Meyn and Tweedie \cite{Meyn1993c}. In particular for Feller processes there have been several attempts to classify their behavior based on the generator or the associated Dirichlet form, see Chapter 6 of Jacob \cite{Jaco2004} and the references given therein. 

In one dimension a transient process either drifts to infinity (i.e. $\lim_{t\to\infty} X_t = +\infty$ or $=-\infty$) or it may be oscillating: $\limsup_{t\to\infty} X_t = + \infty \text{ and } \liminf_{t\to\infty} X_t = -\infty.$\\
An oscillating process may be recurrent, transient or neither of those (cf. Sections \ref{rec-trans-I} and \ref{rec-trans-II} for the definitions). Even for such a simple process as the stable-like process (a Markov process with generator  $-(-\Delta)^{\alpha(x)/2}$ and symbol $|\xi|^{\alpha(x)}$, respectively; see Bass \cite{Bass1988a} for a construction) is the recurrence and transience behavior in general unknown. Besides symmetric $\alpha$-stable L\'evy processes the only processes of this type treated in the literature are processes where $\alpha(\cdot)$ is periodic \cite{Fran2006} or related processes where the generator is a symmetric Dirichlet form  \cite{Uemu2002,Uemu2004}.
%The only treated cases in the literature are: the $\alpha$-stable process (trivial, since the process is just a L\'evy process); the symbol is periodic \cite{Fran2006}; the process is generated by a symmetric Dirichlet form \cite{Uemu2002,Uemu2004}. 
The initial motivation for this paper was to treat the non-symmetric case. But in the following we develop a more general framework.

In Section \ref{rec-trans-I} we introduce a ``local" notion of recurrence and transience for which we will give sufficient conditions in Section \ref{framework}. Afterwards in Section \ref{rec-trans-II} the local notions are linked to the (global) recurrence and transience of the processes. In particular conditions which imply the recurrence-transience dichotomy are given. Furthermore we give a result which allows to compare the behavior of Markov processes which coincide outside some compact ball.
 
The paper closes with an application to stable and stable-like processes.

\section{Recurrence and Transience} \label{rec-trans-I}
We consider time homogeneous strong Markov processes $(\Omega, \Fskript, \Fskript_t, X_t, \theta_t, \Prob_x)$ with c\`adl\`ag paths on $\R^d$ ($d\in\N$), where the filtration $(\Fskript_t)_{t\geq 0}$ satisfies the usual conditions. Note that $(\theta_t)_{t\geq 0}$ is the family of shift operators on $\Omega$, i.e. $X_s(\theta_t(\omega)) = X_{t+s}(\omega)$ for $\omega \in \Omega$.\\
To simplify notation we denote such a process by $(X_t)_{t\geq 0}.$ The state space $\R^d$ will be equipped with the Borel-$\sigma$-algebra $\Bskript(\R^d)$ and sets will be elements of $\Bskript(\R^d)$ if not stated otherwise. For a set $A$ the first entrance time is defined, with the convention $\inf \emptyset = \infty$, by
$$\tau_A:=\inf\{t \geq 0\ |\ X_t \in A\}.$$
Note that $\tau_A$ is a stopping time for any $A\in \Bskript(\R^d)$, since the process is right continuous and adapted, hence progressive. Furthermore for any stopping time $\sigma$ also 
$$\tau_{A,\sigma} := \inf\{ t\geq \sigma\ |\ X_t \in A\}$$
is a stopping time since
$$\{\tau_{A,\sigma} \leq t\} = \bigcup_{s\in\Q \cap [0,t]} \{X_s \in A\}\cap \{\sigma \leq s\} \ \in \Fskript_t$$
(compare \cite{EthiKurt86}, Chapter 2, Prop. 1.5).

Now we define a pointwise (local) notion of recurrence and transience.
\begin{Def} \label{local-defs}
Let $(X_t)_{t\geq 0}$ be $\R^d$-valued process and $b\in\R^d$. With respect to $(X_t)_{t\geq 0}$ the point $b$ is called
\begin{itemize}
\item \textbf{recurrent} if 
$$\Prob_b(\forall T > 0\ \exists t > T:\ X_{t} = b) = 1,$$
\item \textbf{left limit recurrent} if
$$\Prob_b(\forall T > 0\ \exists t > T:\ X_{t-} = b) = 1,$$
\item \textbf{locally recurrent} if 
$$\Prob_b(\liminf_{t\to\infty} |X_t - b| = 0) = 1,$$
\item \textbf{locally transient} if 
$$\Prob_b(\liminf_{t\to\infty} |X_t - b| = 0) < 1,$$
\item \textbf{transient} if 
$$\Prob_b(\liminf_{t\to\infty} |X_t - b| = \infty) = 1.$$
\end{itemize}
\end{Def}
\begin{Rem}
The notion of \textit{local} is meant in a spatial sense, as opposed to a temporal sense.  One would get the latter by transferring the definition of (deterministic) locally recurrent functions (e.g. \cite{Bush1962}) to processes.
\end{Rem}

Note that only for left limit recurrence we need that the paths have left limits, the right continuity is not necessary for these definitions. The reason of introducing left limit recurrence at all, is that our method will not allow to prove recurrence for points but at most left limit recurrence. Nevertheless we have the following Lemma to conclude recurrence for a point.

\begin{Lem} \label{qlc}
Let $(X_t)_{t\geq 0}$ be quasi left continuous, i.e. for every increasing sequence of stopping times $\sigma_n$ with limit $\sigma$:
$$X_{\sigma_n}\xrightarrow{n\to \infty} X_\sigma \ \ \text{ a.s. on }\{\sigma < \infty\}.$$
Then the following implication holds:\\
\centerline{ $b$ is left limit recurrent \ \ $\Rightarrow$\ \  $b$ is recurrent.} 
\end{Lem}
\begin{proof}
Define $\sigma_0:= k \in \N$ and for $n\in\N$
$$\sigma_n:= \inf\{t\geq\sigma_{n-1}\ \big|\ |X_t-b| < \frac{1}{n}\} \ \ \text{ and } \ \  \sigma:= \lim_{n\to\infty} \sigma_n.$$
Clearly $(\sigma_n)_{n\in\N}$ is increasing. Thus $\sigma$ is well defined and 
$$\Prob_b(\sigma < \infty) = 1,$$
since $b$ is left limit recurrent. Note that $\sigma_n$ might be constant for $n$ large, but in this case the process is already in $b.$ In general by the quasi left continuity 
$$\Prob_b(X_\sigma = \lim_{n\to\infty} X_{\sigma_n} = b) = 1$$ 
holds. Since $k$ was arbitrary this yields that $b$ is recurrent.
\end{proof}
 
Further simple consequences of Definition \ref{local-defs} are that (left limit) recurrence implies local recurrence and that we have the dichotomy
\begin{equation}\label{dichotomy}
 b \text{ is either locally recurrent or locally transient}.
\end{equation}

A process $(X_t)_{t\geq 0}$ is point recurrent if and only if all $b\in \R^d$ are recurrent. The other common notions for recurrence and transience of processes do not have such a simple relation to the above local notions. Details will be given in Section \ref{rec-trans-II}.

\section{Overshoots and Markov processes} \label{framework}
In this section we treat for simplicity the case $d=1$, see Remark \ref{extension-dim} for the extension to higher dimensions. Let $(X_t)_{t\geq 0}$ be a process on $\R$ satisfying 
\begin{equation}\label{assume-oscillate}
\limsup_{t\to\infty} X_t = \infty \ \ \text{ and } \ \ \liminf_{t\to\infty} X_t = -\infty \ \ a.s..
\end{equation}
Further assume that there exists some $b\in \R$ such that for the stopping times
$$\tau^b := \inf\{t \geq 0\ |\  X_t \leq b\}\ \ \text{ and } \ \ \sigma^b:= \inf\{t \geq 0\ |\ X_t \geq b\}$$
the process satisfies
\begin{equation}\label{assume-shoot}
\begin{split}
\Prob_x(X_{\tau^b} = b) = 0 &\ \ \ \ \text{ for all } x> b,\\
\Prob_x(X_{\sigma^b} = b) = 0 &\ \ \ \ \text{ for all } x< b,
\end{split}
\end{equation}
i.e. the process almost surely enters $(-\infty,b]$ and $[b,\infty)$ not by hitting $b$. The distributions of $X_{\tau^b}$ and $X_{\sigma^b}$ are called {\bf \textit{overshoot distributions}}. %, also common are the names residual life or excess distribution. 
\begin{Rem}
Note that assumption \eqref{assume-shoot} is not equivalent to assuming that the process is non-creeping. For example consider a compound Poisson process on $\R$ with jump distribution $\frac{1}{2}\delta_{-1}+\frac{1}{2}\delta_{1}$. The process started in 0 is non-creeping but hits $b=1$ with probability one.  
\end{Rem}

Now define $\sigma_0:= 0$ and for each $n\in\N$ set
\begin{align*}
\tau_n&:= \inf \{ t\geq \sigma_{n-1}\ |\ X_t < b\},\\
\sigma_n&:= \inf\{ t\geq \tau_{n}\ |\ X_t > b\}.
\end{align*}
%To keep notation simple we will write $\tau_n$ and $\sigma_n$ for $\tau_n^b$ and $\sigma_n^b$, respectively. %
Note that $\sigma_1$ is always the the first time of passing $b$ from below. Contrary  $\tau_1$ is for the process started in $x>b$ the first time of passing $b$ from above, but $\tau_1 = 0$ for $x<b$.

These stopping times have the following properties.
\begin{Pro} \label{st-properties} Let $x\neq b$, then
\begin{enumerate}[\upshape i)]
\item $\Prob_x(\tau_n < \infty) = 1$ and $\Prob_x(\sigma_n < \infty) = 1$ for all $n\in \N$,
\item $\{X_{\tau_n} < b\} \subset \{\sigma_n > \tau_n\}$,
\item $\Prob_x(X_{\tau_n} < b ) = 1$ implies $\Prob_x(X_{\sigma_n} > b) = 1$,  
\item $\Prob_x(X_{\sigma_{n}}>b,X_{\tau_n} < b,\ \forall n\in\N) = 1$,
\item $\Prob_x(\sigma_{n-1} < \tau_n < \sigma_n,\ \forall n\in\N) = 1$.
\end{enumerate}
\end{Pro}
\begin{proof}
\begin{enumerate}[i)]
\item By \eqref{assume-oscillate} the process will pass $b$ infinitely often almost surely, i.e. $\tau_n$ and $\sigma_n$ are finite almost surely.
\item Let $\omega \in \{X_{\tau_n} < b\}$. Then by the right continuity there exists an $\varepsilon_\omega > 0$ such that $X_{\tau_n + \varepsilon_\omega}(\omega) < b$, since $(X_t)_{t\geq 0}$ is c\`adl\`ag. Thus $\sigma_n(\omega) \geq \tau_n(\omega) + \varepsilon_\omega$, i.e. 
$$\sigma_n(\omega) > \tau_n(\omega).$$
\item First note that $\Prob_x(X_{\tau_n} < b ) = 1$ implies by ii) that $\Prob_x(\sigma_n > \tau_n) = 1,$ and $\tau_n$ is a finite stopping time by i). By the right continuity $\{X_{\sigma_1} = b\}$ contains all paths which enter $(b,\infty)$ continuously from $b$ and $\{X_{\sigma^b} = b\}$ contains all paths which enter $[b,\infty)$ at $b.$ Thus $\{X_{\sigma_1}=b\} \subset \{X_{\sigma^b} = b\}$, i.e.
$$\Prob_y(X_{\sigma_1}=b) \leq \Prob_y(X_{\sigma^b}=b) = 0,$$
which implies $\Prob_y(X_{\sigma_1}>b) = 1.$
Now for $y<b$ the strong Markov property (note: $\sigma_n = \sigma_1 \circ \theta_{\tau_n}$) yields
$$\Prob_x(X_{\sigma_n}>b|X_{\tau_n} = y) = \Prob_y (X_{\sigma_1} > b) =  1.$$
Then 
\begin{equation*}
\begin{split}
\Prob_x(X_{\sigma_n} > b) & = \int_{(-\infty, b]} \Prob_x(X_{\sigma_n} > b| X_{\tau_n} = y) \ \Prob_x(X_{\tau_n} \in dy)\\
& = \int_{(-\infty,b]} 1 \ \Prob_x(X_{\tau_n} \in dy) = 1.
\end{split}
\end{equation*}
\item Analogously to ii) and iii) one gets:
\begin{itemize}
\item[ii*)] $\{X_{\sigma_n} > b\} \subset \{ \tau_{n+1} > \sigma_n\},$
\item[iii*)] $\Prob_x(X_{\sigma_n}>b) = 1$ implies $\Prob_x(X_{\tau_{n+1}} < b) = 1,$
\end{itemize}
and further
$$\Prob_x(X_{\tau_1}< b) = 1$$
holds. Thus repeated applications of iii) resp. iii*) yield 
$$\Prob_x(X_{\tau_n}<b) = \Prob_x(X_{\sigma_n} > b) = 1 \ \text{ for each }n\in \N.$$
Thus 
$$ \Prob_x(X_{\tau_n}< b, X_{\sigma_n}>b,\ \forall n\in \N) = 1$$
as a countable intersection of sets of measure one.
\item This is a consequence of ii), ii*) and iv). \qedhere
\end{enumerate}
\end{proof}

Now define for $x>b$ on the set $\{\sigma_{n-1} < \tau_n < \sigma_n,\ \forall n\in\N\}$, which has probability one by Proposition \ref{st-properties} v), the sequence $(Y_n)_{n\geq 0}$ by
$$Y_n:=X_{\sigma_n}$$
and note that by the strong Markov property for $B\in \Borel(\R)$
$$\Prob(Y_n\in B\ |\ Y_{n-1} = x)=\Prob(X_{\sigma_n}\in B\ |\ X_{\sigma_{n-1}} = x) = \Prob(X_{\sigma_1}\in B\ |\ X_0 = x)=\Prob(Y_1\in B\ |\ Y_0 =x),$$
i.e. $(Y_n)_{n\geq 0}$ is a Markov chain on $(b,\infty)$. This Markov chain captures only the \textit{first set of countably many overshoots} passing $b$ from $(-\infty,b)$ of the process $(X_t)_{t\geq 0}$, since the times $(\sigma_n)_{n\geq 0}$ are strictly increasing but possibly bounded.  

Nevertheless this Markov chain can be used to determine the local recurrence/transience behavior of $(X_t)_{t\geq 0}$ by the following theorem. 

\begin{Thm} \label{recurrence-transience-thm} Let $(X_t)_{t\geq 0}$ and $(Y_n)_{n\geq 0}$ be as defined above.
\begin{enumerate}[\upshape i)]
\item If $\Prob_x(\lim_{n\to\infty} Y_n = \infty)=1$ for all $x>b$ and there exists $r,R >0$ and $c<1$ such that
\begin{equation}
\label{local-transience-estimate} \sup_{\substack{y\in[b-r,b+r]\\y\neq b}} \Prob_y(X_{\sigma_1} > b+R)< c
\end{equation}
then $b$ is locally transient.
\item If $\Prob_x(\liminf_{n\to\infty} Y_n = b)=1$  for all $x>b$ then $b$ is locally recurrent.
\item If $\Prob_x(\lim_{n\to\infty} Y_n = b)=1$  for all $x>b$ and there exists $r',R' >0$ and $c<1$ such that
\begin{equation}
\label{leftlimit-recurrence-estimate} \sup_{y\geq b+r'} \Prob_y(X_{\sigma_1} < b+R')< c
\end{equation}
then $b$ is left limit recurrent.
\end{enumerate}
\end{Thm}
\begin{Rem} Roughly speaking, condition \eqref{local-transience-estimate} ensures that the overshoots \textit{represent} the whole process, whereas condition \eqref{leftlimit-recurrence-estimate} ensures that the limit $b$ is reached in finite time.  
The following two examples show these conditions cannot be removed.
\begin{enumerate}
\item Let $(N_t)_{t\geq 0}$ be a Poisson process and $(\tilde{X}_n)_{n\geq 0}$ be a Markov chain with transition distribution 
$$\Prob(\tilde{X}_1 \in dy\ |\ \tilde{X}_0 = x) = \begin{cases} \delta_\frac{1}{x}(dy)\ &\text{ for }|x|>1,\\
\delta_{-\frac{1 + |x|}{x}}(dy) &\text{ for } 0< |x|\leq 1,\\
\delta_1(dy) &\text{ for } x = 0.
\end{cases}$$
The Markov chain is in fact deterministic and, when started in 0, the chain moves as
$$0,1,-2,-\frac{1}{2},3,\frac{1}{3},-4,-\frac{1}{4},\ldots .$$
Now the chain subordinated by the Poisson process is a c\'adl\'ag time homogeneous strong Markov process satisfying \eqref{assume-oscillate} and \eqref{assume-shoot} for $b=0.$ Furthermore 0 is locally recurrent and thus not locally transient. The associated chain of overshoots is deterministic, especially for $x\in(0,1]:$
$$Y_0 = x,\ Y_1=\frac{1}{x}+2 \text{ and for $n\in\N$ } Y_n= \frac{1}{x}+2n,$$ 
i.e. $\lim_{n\to\infty} Y_n = \infty$ and 
$$\forall R,r > 0:\ \sup_{y\in[-r,r],y\neq 0} \Prob_y(X_{\sigma_1}>R) \geq \sup_{y\in(0,r]}\Prob_y(Y_1 >R) = 1.$$

\item Changing the transition distribution to
$$\Prob(\tilde{X}_1 \in dy\ |\ \tilde{X}_0 = x) = \begin{cases} \delta_{-\frac{1}{x}}(dy)\ &\text{ for }|x|>1,\\
\delta_{\frac{1+|x|}{x}}(dy) &\text{ for } 0< |x|\leq 1,\\
\delta_1(dy) &\text{ for } x = 0.
\end{cases}$$
yields that the chain started in 0 moves as
$$0,1,2,-\frac{1}{2},-3,\frac{1}{3},4,-\frac{1}{4},-5,\ldots .$$
Thus for the chain subordinated by the Poisson process 0 is locally recurrent but not left limit recurrent (in finite time).
For the associated jump chain for $x>1$ we find
$$Y_0 = x,\ Y_1=\frac{1}{x+1} \text{ and in general } Y_n= \frac{1}{x+2n-1},$$ 
i.e. $\lim_{n\to\infty} Y_n = 0$ and $\forall R,r > 0:\ \sup_{y\geq r}\Prob_y(X_{\sigma_1} < R) = \sup_{y\geq r}\Prob_y(Y_1 < R) = 1.$
\end{enumerate}
\end{Rem}
\begin{proof}[Proof of Theorem \ref{recurrence-transience-thm}]
\begin{enumerate}[i)]
\item By \eqref{assume-oscillate} $(X_t)_{t\geq 0}$ does not explode in finite time. This and $\infty = \lim_{n\to\infty} Y_n = \lim_{n\to\infty} X_{\sigma_n}$ a.s. imply that $\sigma_n \to \infty$ almost surely. Let $r, R$ and $c$ be as in \eqref{local-transience-estimate}. Now fix $\varepsilon > 0.$ Then there exists a $N>0$ such that 
$$\forall n\geq N:\ \Prob_x(X_{\sigma_n} > R+b) \geq 1-\varepsilon,$$
since $\lim_{n\to\infty} X_{\sigma_n} = \infty.$\\
 Let $n\geq N$ and define $\nu_n$ as the time of the first visit to $B:=[b-r,b+r]\backslash\{b\}$ after time $\sigma_n,$ i.e.
$$\nu_n:= \inf\{t\geq \sigma_n\ |\ X_t \in B\}$$
and $\sigma_k$ be the time of the first jump into $(b,\infty)$ after $\nu_n$, i.e.
$$k:= \inf\{l \in \N\ |\ \sigma_l > \nu_n\}.$$

Now suppose $b$ is locally recurrent. An overshoot hits $b$ with probability zero, thus the local recurrence of $b$ implies that that $\Prob_x(\nu_n < \infty) = 1$ and $\Prob_x(X_{\nu_n} \in B) = 1.$  Thus $\Prob_x(k<\infty) = 1$ and $\sigma_k = \sigma_1 \circ \theta_{\nu_n}$, where $\theta_{\nu_n}$ is the shift operator corresponding to $\nu_n$. Then the strong Markov property yields
\begin{equation*}
\begin{split}
1-\varepsilon &\leq \Prob_x(X_{\sigma_k} > R+b) \\
&=\int_{B} \Prob_x(X_{\sigma_k}>R+b| X_{\nu_n} = y)\ \Prob_x(X_{\nu_n}\in dy) \\
&= \int_{B} \Prob_y(X_{\sigma_1} > R+b) \ \Prob_x(X_{\nu_n}\in dy) \\
&\leq \sup_{y\in B} \Prob_y(X_{\sigma_1} > R+b)< c  <1,
\end{split}
\end{equation*}
which is a contradiction. Thus $b$ is locally transient.

\item Let $\liminf_{n\to\infty} Y_n = b$ almost surely. If $\sigma_n \to \infty$ a.s. the statement is obvious. In general let $\varepsilon > 0,$ $T> 0$ and 
$$\eta_T := \inf\{t\geq T\ |\ X_t \in [1,\infty)\}.$$
By \eqref{assume-oscillate} for all $y\in\R$ we have $\Prob_y(\eta_T < \infty) = 1.$ Thus for $x>0$ the strong Markov property yields
\begin{equation*}
\begin{split}
\Prob_x(\exists t>T \ :\ |X_t-b| <\varepsilon) &\geq \Prob_x(\exists t>0 \ :\ |X_{t+\eta_T}-b| <\varepsilon)\\
&= \int_{[1,\infty)} \Prob_x( \exists t> 0\ :\ |X_{t+\eta_T} - b| < \varepsilon\ \big|\  X_{\eta_T} = y)\ \Prob_x(X_{\eta_T} \in dy)\\
&= \int_{[1,\infty)} \Prob_y( \exists t> 0\ :\ |X_{t} - b| < \varepsilon)\ \Prob_x(X_{\eta_T} \in dy)\\
&\geq \int_{[1,\infty)} \Prob_y( \exists n \in \N \ :\ |Y_{n} - b| < \varepsilon)\ \Prob_x(X_{\eta_T} \in dy) = 1.
\end{split}
\end{equation*}
Since $T$ and $\varepsilon$ where arbitrary this implies that $b$ is locally recurrent.
\item If $(\sigma_n)_{n\in\N}$ is a.s. bounded then $b$ is reached as the left limit at least once and the same argument as in part $ii)$ implies that $b$ is left limit recurrent.\\ Otherwise set $\sigma_\infty := \lim_{n\to\infty} \sigma_n$ and let $r',R'$ and $c$ be as in \eqref{leftlimit-recurrence-estimate}. Further let $\varepsilon > 0$ and $N>0$ such that for all $n\geq N$
$$\Prob_x(X_{\sigma_n}<b+R') \geq 1-\varepsilon,$$
such an $N$ exists since $\lim_{n\to\infty}Y_n = b$ a.s.. Now let $n\geq N$ and define  $\nu_n$ as the time of the first visit to $(b+r',\infty)$ after time $\sigma_n,$ i.e.
$$\nu_n:= \inf\{t\geq \sigma_n| X_t \geq b + r'\}$$
and 
$$k:= \inf\{l\in\N\ |\ \sigma_l > \nu_n\}.$$
Note that $\sigma_\infty = \infty$ with positive probability but in general not almost surely. Thus only on $\{\sigma_k > \nu_n\}$ the stopping time $\sigma_k$ is the time of the first jump into $(b,\infty)$ after $\nu_n$, i.e. on this set 
$\sigma_k = \sigma_1 \circ \theta_{\nu_n}$
holds.
Now $\One_{\{\sigma_k>\nu_n\}}$ is $\Fskript_{\nu_n}$ measurable and the strong Markov property by conditioning on $\Fskript_{\nu_n}$ (the $\sigma$-algebra associated with $\nu_n$) yields
\begin{equation*}
\begin{split}
1-\varepsilon &\leq \Prob_x(X_{\sigma_k} < b+R') \\
&=\Prob_x(X_{\sigma_k} < b+R',\ \sigma_k \leq\nu_n) + \Prob_x(X_{\sigma_k} < b +R',\ \sigma_k>\nu_n)\\
&\leq \Prob_x(\sigma_k \leq \nu_n) + \E_x(\E_x(\One_{\{X_{\sigma_1}(\theta_{\nu_n}) < b+R'\}}\One_{\{\sigma_k > \nu_n\}}\ |\ \Fskript_{\nu_n}))\\
&= \Prob_x(\sigma_k \leq \nu_n) + \E_x(\One_{\{\sigma_k > \nu_n\}} \Prob_{X_{\nu_n}}(X_{\sigma_1} < b+R'))\\
&< \Prob_x(\sigma_k \leq \nu_n) + c\Prob(\sigma_k>\nu_n) \leq \tilde{c} < 1,
\end{split}
\end{equation*}
which is a contradiction, since $\varepsilon$ was arbitrary. Thus $(\sigma_n)_{n\in\N}$ is bounded. \qedhere
\end{enumerate}
\end{proof}

\begin{Rem} \label{extension-dim}
In order to use this approach for $d>1$ one has to replace $(-\infty,b]$ and $[b,\infty)$ by parts of the state space separated by a $(d-1)$-dimensional hyperplane. Furthermore \eqref{assume-oscillate} has to be reformulated, such that it ensures that the process passes the hyperplane infinitely often and reaches an arbitrary large distance to the hyperplane. Then analogous to \eqref{assume-shoot} it has to be required that the up/down shoots with respect to the hyperplane do not hit it. With this an analogue to Proposition \ref{st-properties} holds. Also an analogous result to Theorem  \ref{recurrence-transience-thm} can be proved. For part i) condition \eqref{local-transience-estimate} has to be defined with respect to the hyperplane and the limit of the distance of the overshoots to the hyperplane should become arbitrary large with probability 1, part ii) for $b\in\R^d$ is analogous to the one dimensional case and part iii) requires again a reformulation of \eqref{leftlimit-recurrence-estimate} in terms of the hyperplane. 

But note that for $d>1$ the set of cases where the theorem does not lead a conclusion will be considerably larger than in one dimension, since the transience part only considers deviations which are (in a sense) orthogonal to the hyperplane.
\end{Rem}

\section{Recurrence and Transience of Processes}
\label{rec-trans-II}
In this section we will link local recurrence and local transience to the notion of recurrence and transience for processes, as used by Meyn and Tweedie e.g. in \cite{Meyn1993c} (our presentation is partly motivated by \cite{Stra1997}). Note that all results of this section would also hold if we weaken our assumption on the processes from \textit{c\`adl\`ag} to only \textit{right continuous}.

By $\lambda$ we denote the Lebesgue measure.
\begin{Def}
A process $(X_t)_{t\geq 0}$ on $\R^d$ is called
\begin{itemize}
\item \textbf{$\mathbf{\lambda}$-irreducible} if
$$\lambda(A) > 0 \Rightarrow \E_x\left(\int_0^\infty \One_A(X_t) \ dt\right) > 0 \text{ for all } x,$$

\item \textbf{recurrent} with respect to $\lambda$ if 
$$\lambda(A) > 0 \Rightarrow \E_x\left(\int_0^\infty \One_A(X_t) \ dt\right) = \infty \text{ for all } x,$$
\item \textbf{Harris recurrent} with respect to $\lambda$ if 
$$\lambda(A) > 0 \Rightarrow \Prob_x\left(\int_0^\infty \One_A(X_t) \ dt =\infty\right) = 1 \text{ for all } x,$$
\item \textbf{transient} if there exists a countable cover of $\R^d$ with sets $A_j$ such that for each $j$ there is a finite constant $M_j>0$ such that:
$$\E_x\left(\int_0^\infty \One_{A_j}(X_t) \ dt\right) < M_j,$$
\item a \textbf{$\mathbf{T}$-model} if for some probability measure $\mu$ on $[0,\infty)$ there exists a kernel $T(x,A)$ with $T(x,\R^d)>0$ for all $x$ such that the function $x\mapsto T(x,A)$ is lower semi-continuous for all $A\in\Bskript(\R^d)$ and
$$ \int_0^\infty \E_x(\One_A(X_t)) \ \mu(dt) \geq T(x,A)$$
holds for all $x$, $A\in \Bskript(R^d)$.
\end{itemize}
\end{Def}

We start with the recurrence-transience-dichotomy for $\lambda$-irreducible $T$-models.
\begin{Thm} \label{classical-dichotomy}
Let $(X_t)_{t\geq 0}$ be a $\lambda$-irreducible $T$-model, then it is either Harris recurrent or transient.
\end{Thm}
\begin{proof}
(Compare with the proof of Prop. 3.1 in \cite{Stra1997}.) A $\lambda$-irreducible process is by Thm. 2.3 in \cite{Twee1994} either recurrent or transient. In the case of recurrence the reference measure is the so called maximal irreducible measure, but this yields in our case especially recurrence with respect to $\lambda.$\\
Now suppose the process is recurrent with respect to $\lambda$ then for all $x\in\R^d$ and all $\varepsilon >0$
$$\E_x\left(\int_0^\infty \One_{B_\varepsilon(x)}(X_t) \ dt\right) = \infty \ \ \text{ where } B_\varepsilon(x) = \{y\in\R^d \ | \ |x-y|< \varepsilon\}$$
i.e. each $x\in\R^d$ is \textit{topological recurrent} (cf. Sec. 4 \cite{Twee1994})). Thus by Thm. 4.2 in \cite{Twee1994} the whole state space $\R^d$ is a maximal Harris set, that means there exists a measure $\phi$ on $\Bskript(\R^d)$ such that $(X_t)_{t\geq 0}$ is Harris recurrent with respect to $\phi.$ Now $\phi = \mu R$ (cf. proof of Thm. in 2.4 \cite{Meyn1993b} and the proof of Prop. 3.1 in \cite{Stra1997}) for some non trivial measure $\mu$ and a kernel $R$ which satisfies
$$\lambda(A)>0 \ \Rightarrow \ R(x,A) >0\ \ \forall x\in\R^d.$$ 
Thus $(X_t)_{t\geq 0}$ is Harris recurrent with respect to $\lambda.$
\end{proof}

Now we can state the main theorem of this section which links the local notions introduced in Section \ref{rec-trans-I} to the stability of the process.

\begin{Thm} \label{dichotomy-tweedie}
Let $(X_t)_{t\geq 0}$ on $\R^d$ be a $\lambda$-irreducible $T$-model, then
\begin{enumerate}[\upshape i)]
\item $\exists$ $b$ which is locally recurrent $\Leftrightarrow$ $(X_t)_{t\geq 0}$ is Harris recurrent.
\item $\exists$ $b$ which is locally transient $\Leftrightarrow$ $(X_t)_{t\geq 0}$ is transient.
\end{enumerate}
\end{Thm}
\begin{proof}
By Theorem \ref{classical-dichotomy} the process is either Harris recurrent or transient.   Thus it is enough to prove the equivalence in i) since also local recurrence and local transience are complementary.

As in the previous proof, a point $x$ is called topologically recurrent if $\E_x(\int_0^\infty \One_A(X_t)\ dt) = \infty$ for all neighborhoods $A$ of $x$. Note that for a $\lambda$-irreducible process each point is reachable, i.e. for every $x$ and every neighborhood $A$ we have $\Prob_x(\tau_A <\infty) >0$. Thus Thm. 4.1 in \cite{Twee1994} yields for a $\lambda$-irreducible T-model:
$$\exists\ b\ \text{topologically recurrent} \Leftrightarrow (X_t)_{t\geq 0} \text{ is recurrent }.$$
   
Now assume $b$ is locally recurrent. For any neighborhood $A$ of $b$ we find a open ball with center $b$ and radius $\varepsilon>0$ such that $B_\varepsilon(b) \subset A$. The local recurrence implies that the process hits $B_\frac{\varepsilon}{2}(b)$ with probability one, also after arbitrary large times, i.e. for all $R>0$
$$\Prob_b(\exists t>R:\ X_t \in B_\frac{\varepsilon}{2}(b))=1.$$
Furthermore since $X_t$ is right continuous the average time spent in $B_\varepsilon(b)$ after hitting $B_\frac{\varepsilon}{2}(b)$ is positive, i.e.
$$0 < \inf_{y\in B_\frac{\varepsilon}{2}(b)} \E_y(\tau_{\R^d\backslash B_\varepsilon(b)}).$$
Thus we get
$$\E_b\left(\int_0^\infty \One_A(X_t)\ dt\right) \geq \E_b\left(\int_0^\infty \One_{B_\varepsilon}(X_t)\ dt\right) \geq \infty,$$
i.e. $b$ is topological recurrent. Therefore $(X_t)_{t\geq 0}$ is recurrent. By the dichotomy we get that in fact $(X_t)_{t\geq 0}$ is Harris recurrent, since it is not transient. 

On the other Hand, let $(X_t)_{t\geq 0}$ be Harris recurrent. Thus
$$\Prob_x\left(\int_0^\infty \One_A(X_t) \ dt =\infty\right) = 1 \text{ for all } x \text{ and all } A \text{ with } \lambda(A)>0$$
holds and especially the path returns into $B_\varepsilon(b)$ for any $\varepsilon>0$ after any time, i.e. $b$ is locally recurrent.   
\end{proof}

We further recall the following theorem, which provides some way to check that $(X_t)_{t\geq 0}$ is a T-model. 

\begin{Thm}[Thm. 5.1 and Thm. 7.1 in \cite{Twee1994}] \label{thm-T-model}
\begin{enumerate}[\upshape i)]
\item $(X_t)_{t\geq 0}$ is a T-model, if every compact set $C$ is petite, i.e. there exists a probability measure $\mu$ on $[0,\infty)$ and a non-trivial measure $\nu$ on $\R^d$ such that
$$ \int_0^\infty \E_x(\One_A(X_t)) \ \mu(dt) \geq \nu(A) \ \text{ for all } x\in C \text{ and all } A.$$
\item Let $(X_t)_{t\geq 0}$ be $\lambda$-irreducible and $x \mapsto \E_x(f(X_t))$ be continuous for all  continuous and bounded functions $f$, then $(X_t)_{t\geq 0}$ is a T-model.
\end{enumerate}
\end{Thm}

Part ii) in particular shows that every $\lambda$-irreducible $C_b$-Feller process is a T-model, and note that \cite{Schi1998} gives necessary and sufficient conditions for a $C_\infty$-Feller process to be also $C_b$-Feller. 

Useful for applications is the following theorem which gives sufficient criteria for a process to be a $\lambda$-irreducible T-model.

\begin{Thm}
Let $(X_t)_{t\geq 0}$ be a process on $\R^d$ and denote its transition probabilities by
$$P_t(x,A) := \Prob_x(X_t\in A).$$
Then
\begin{enumerate}[\upshape i)]
\item $(X_t)_{t\geq 0}$ is $\lambda$-irreducible if
\begin{equation}
\label{lebesgue-transition}
\lambda(A) > 0 \Rightarrow P_t(x,A) > 0 \text{ for all } t>0,x\in\R^d,
\end{equation}
\item $(X_t)_{t\geq 0}$ is a $\lambda$-irreducible T-model if \eqref{lebesgue-transition} holds and there exits a compact set $K\subset[0,\infty]$ and a non trivial measure $\nu$ such that for all compact sets $C\subset \R^d$ 
\begin{equation}
\label{lebesgue-T}
\inf_{t\in K} \inf_{x\in C} P_t(x,A) \geq \nu(A) \text{ for all } A \in \Bskript(\R^d).
\end{equation}
\end{enumerate}
Further, a special case of {\upshape ii)}:
\begin{enumerate}[\upshape i)]
\item[\upshape iii)] $(X_t)_{t\geq 0}$ is a $\lambda$-irreducible T-model if the transition probability $P_t(x,.)$ is the sum of a, possibly trivial, discrete measure and a measure which has a (sub-)probability density $\tilde{p}_t(x,y)$ with respect to $\lambda$ such that
\begin{eqnarray} \label{density-condition-a}
\tilde{p}_t(x,y) > 0& \text{ for all } x,y\in\R^d,\ t>0,\\
\label{density-condition-b}
 \inf_{t\in[1,2]}\inf_{x\in C} \tilde{p}_t(x,y) >0& \text{ for all } y\in R^d \text{ and all compact sets } C. 
\end{eqnarray}
\end{enumerate}
\end{Thm}
\begin{proof}
Assume \eqref{lebesgue-transition} holds and let $A$ be such that $\lambda(A)>0$. Then
$$\Prob_x(\tau_A < \infty) \geq P_t(x,A) > 0 \text{ for any }t>0$$
and thus by Prop. 2.1 of \cite{Meyn1993c} the process $(X_t)_{t\geq 0}$ is $\phi$-irreducible with
$$\phi(.) := \lambda(.) \int_{[0,\infty)} e^{-t} P_t(x,.)\ dt.$$
But clearly for $A$ with $\lambda(A)>0$ also 
$$\int_{[0,\infty)} e^{-t} P_t(x,A)\ dt > 0$$
holds. Therefore $\phi$ is equivalent to $\lambda$, i.e. $(X_t)_{t\geq 0}$ is $\lambda$-irreducible.

If further \eqref{lebesgue-T} holds then Theorem \ref{thm-T-model} part i) with $\mu(dt) = e^{-t}dt$ implies that $(X_t)_{t\geq 0}$ is a T-model.

For part iii) note that \eqref{density-condition-a} implies that \eqref{lebesgue-transition} holds and \eqref{density-condition-b} implies that \eqref{lebesgue-T} holds with $\nu$ being a subprobability measure with density $\inf_{t\in[1,2]}\inf_{x\in C} \tilde{p}_t(x,.) e^{-2}$.
\end{proof}

We give a further characterization of recurrence and transience in this context, which shows that it is in fact enough to know the behavior of the process outside some compact set.

\begin{Thm} \label{recurr-trans-stoppingtime}
Let $(X_t)_{t\geq 0}$ be $\lambda$-irreducible T-model, $R$ be some positive constant and $\overline{B_R(0)}$ denote the closed ball centered at 0 with radius $R$, then
\begin{enumerate}[\upshape i)]
\item $\forall x: \Prob_x\left(\tau_{_{\overline{B_R(0)}}} < \infty\right) = 1 \iff (X_t)_{t\geq 0}$ is Harris recurrent.
\item $\exists x: \Prob_x\left(\tau_{_{\overline{B_R(0)}}} < \infty\right) < 1 \iff (X_t)_{t\geq 0}$ is transient.
\end{enumerate}
\end{Thm}
\begin{proof}
Given a $\lambda$-irreducible T-model then by Thm. 5.1 in \cite{Twee1994} every compact set is petite. Thus Thm. 3.3 in \cite{Meyn1993b} implies  ``$\Rightarrow$" of i).

For ii) ``$\Rightarrow$" note that $\lambda(B_R(0)) > 0$. Thus $(X_t)_{t\geq 0}$ cannot be Harris recurrent and the dichotomy implies that it is transient.

Harris recurrence and transience are complementary and so are the left hand sides of i) and ii). Thus the ``$\Leftarrow$" directions hold.   
\end{proof}

In fact the Theorem \ref{recurr-trans-stoppingtime} shows that processes which coincide outside a ball have the same recurrence and transience behavior, respectively.

\begin{Cor} \label{tail-comparison}
Let $(X_t)_{t\geq 0}$ and $(Y_t)_{t\geq 0}$ be $\lambda$-irreducible T-models. If there exists an $R>0$ such that
$$\tau^{X}_{_{\overline{B_R(0)}}} \stackrel{d}{=} \tau^{Y}_{_{\overline{B_R(0)}}} \text{ for all } X_0 = Y_0 = x \in \R^d\backslash\overline{B_R(0)}$$
then $(X_t)_{t\geq 0}$ and $(Y_t)_{t\geq 0}$ have the same recurrence/transience behavior.\\
Here $\tau^X$ and $\tau^Y$ are the entrance times corresponding to $X_t$ and $Y_t$, respectively and $\stackrel{d}{=}$ denotes equality in distribution.
\end{Cor}
\begin{proof}
In the setting of Theorem \ref{recurr-trans-stoppingtime} we find
$$\Prob_x\left(\tau_{_{\overline{B_R(0)}}} = 0\right) = 1 \text{ for all } x\in \overline{B_R(0)}.$$
This shows that Theorem \ref{recurr-trans-stoppingtime} ii) might only hold for some $x \in \R^d\backslash\overline{B_R(0)}$, i.e. only the distributions of $\tau_{\overline{B_R(0)}}\ $ for $x \in \R^d\backslash\overline{B_R(0)}$ need to be checked. Thus, if these distributions coincide for two processes, Theorem \ref{recurr-trans-stoppingtime} yields the same behavior.
\end{proof}

\section{$\alpha$-stable and stable-like Processes}\label{stable-like}
\newcommand{\const}[1]{\frac{\sin(\frac{#1 \pi}{2})}{\pi}}
Let $(X_t)_{t\geq 0}$ be a real valued symmetric $\alpha$-stable process, i.e. it is a L\'evy process with characteristic exponent $|\xi|^\alpha$ with $\alpha \in (0,2)$. In particular it is a time homogeneous strong Markov process with c\`adl\`ag paths. Note that $(X_t)_{t\geq 0}$ sampled at integer times $(X_n)_{n\in N_0}$ is a symmetric random walk and \eqref{assume-oscillate} holds. Define $\sigma^b$ and $\tau^b$ as in Section \ref{framework}, i.e.
$$\tau^b := \inf\{t \geq 0\ |\  X_t \leq b\}\ \ \text{ and } \ \ \sigma^b:= \inf\{t \geq 0\ |\ X_t \geq b\}.$$
In 1958 Ray \cite{Ray1958} showed that for $b>0$ 
$$\Prob_0(X_{\sigma^{b}} \in dy ) = \const{\alpha} \frac{1}{y} \left(\frac{b}{y-b}\right)^{\frac{\alpha}{2}} \One_{[b,\infty)}(y)\ dy
$$
and in particular for $0<\alpha<2$
$$\Prob_0(X_{\sigma^{b}} = b) = 0.$$
The translation invariance of $(X_t)_{t\geq 0}$ yields for all $b$ 
\begin{equation}
\label{alpha-upshoot}
\Prob_x(X_{\sigma^{b}} \in dy ) = \Prob_0(X_{\sigma^{b-x}}+x \in dy ) = \frac{\sin(\frac{\alpha \pi}{2})}{\pi} \frac{1}{y-x} \left(\frac{b-x}{y-b}\right)^{\frac{\alpha}{2}} \One_{[b,\infty)}(y)\ dy \text{ for } x<b
\end{equation}
and the symmetry yields 
\begin{equation}
\label{alpha-downshoot}
\Prob_x(X_{\tau^{b}} \in dy ) = \Prob_{-x}(-X_{\sigma^{-b}} \in dy ) = \frac{\sin(\frac{\alpha \pi}{2})}{\pi} \frac{1}{x-y} \left(\frac{x-b}{b-y}\right)^{\frac{\alpha}{2}} \One_{(-\infty,b]}(y)\ dy \text{ for } x>b.
\end{equation}
In particular \eqref{assume-shoot} is satisfied.

Note that by the translation invariance for any $b$
$$\text{for }x<0:\ \Prob_x(X_{\sigma^0}<r) = \Prob_{x+b}(X_{\sigma^b}< r+b),$$
$$\text{for }x>0:\ \Prob_x(X_{\tau^0}<r) = \Prob_{x+b}(X_{\tau^b}< r+b).$$
Thus we will for simplicity only consider the case $b=0$ in the sequel and define the upwards-overshoot density $u$ and the downwards-overshoot density $v$ for $\alpha \in (0,2)$ by
$$\text{for } x<0: u_\alpha(x,y) := \frac{\sin(\frac{\alpha \pi}{2})}{\pi} \frac{1}{y-x} \left(-\frac{x}{y}\right)^\frac{\alpha}{2} \One_{[0,\infty)}(y)$$
$$\text{for } x>0: v_\alpha(x,y) :=  \frac{\sin(\frac{\alpha \pi}{2})}{\pi} \frac{1}{x-y} \left(-\frac{x}{y}\right)^\frac{\alpha}{2} \One_{(-\infty,0]}(y)$$
We will write $X \sim f$ for a random variable $X$ with density $f$.

\begin{Lem} \label{stable-overshoots-lem} Let $\alpha$, $\beta \in (0,2)$ and $U \sim u_\alpha(-1,\cdot)$ and $V \sim v_\beta(1,\cdot)$ be independent. Then 
\begin{enumerate}[\upshape i)]
\item the overshoot densities satisfy for $y\in\R$
$$\text{for }x<0:\ u_\alpha(x,y) = -\frac{1}{x} u_\alpha(-1,-\frac{y}{x})$$
and
$$\text{for }x>0:\ v_\beta(x,y) = \frac{1}{x} v_\beta(1,\frac{y}{x}),$$
\item for arbitrary probability densities $f$ on $[0,\infty)$ and $g$ on $(-\infty,0]$, and random variables $F \sim f$, $ G\sim g$ independent of $V$ and $U$ respectively, it holds that (for $s\in\R$)
$$\Prob(F V \leq s) = \int_{-\infty}^s \int_{-\infty}^\infty f(x) v_\beta(x,y)\ dx \ dy$$
and
$$\Prob(- G U \leq s) = \int_{-\infty}^s \int_{-\infty}^\infty g(x) u_\alpha(x,y)\ dx\ dy,$$
\item for $r\in\R$
$$\E(U^r) = \begin{cases}\displaystyle \frac{\sin\left(\frac{\alpha \pi}{2}\right)}{\sin\left(\frac{(\alpha-2r)\pi}{2}\right)} &\text{ for } \frac{\alpha}{2} - 1 < r < \frac{\alpha}{2}, \\
\infty & \text{ otherwise},  \end{cases}
$$
and
$$\E((-VU)^r) = \begin{cases}\displaystyle \frac{\sin\left(\frac{\alpha \pi}{2}\right)\sin\left(\frac{\beta \pi}{2}\right)}{\sin\left(\frac{(\alpha - 2r) \pi}{2}\right) \sin\left(\frac{(\beta-2r) \pi}{2}\right)} & \text{ for } \frac{\alpha\lor\beta}{2}-1 < r < \frac{\alpha\land\beta}{2},\\
\infty & \text{ otherwise,}\end{cases}$$
\item for $\alpha+\beta \neq 2$ there exists a moment of a downwards-overshoot followed by an upwards-overshoot which is less than 1, i.e.
$$\alpha + \beta < 2:\ \exists r<0:\ \E(-VU)^r < 1,$$
$$\alpha + \beta > 2:\ \exists r>0:\ \E(-VU)^r < 1,$$
and for $\alpha+\beta = 2$ there is a symmetry:
$$\forall s:\ \Prob(-VU\leq s) = \Prob((-VU)^{-1} \leq s).$$

\end{enumerate}
\end{Lem}
\begin{proof}
\begin{enumerate}[i)]
\item For $x<0$ 
\begin{equation*}
\begin{split}
-\frac{1}{x} u_\alpha\left(-1,-\frac{y}{x}\right) &= \frac{\sin(\frac{\alpha \pi}{2})}{\pi} \frac{1}{-x \left(- \frac{y}{x} +1\right)} \left(\frac{1}{- \frac{y}{x}}\right)^\frac{\alpha}{2} \One_{[0,\infty)}\left(-\frac{y}{x}\right)\\
 &= \frac{\sin(\frac{\alpha \pi}{2})}{\pi} \frac{1}{y-x} \left(-\frac{x}{y}\right)^\frac{\alpha}{2} \One_{[0,\infty)}(y)\\
 & = u_\alpha(x,y)
\end{split}
\end{equation*}
holds and analogously for $x>0$
\begin{equation*}
\begin{split}
\frac{1}{x} v_\beta\left(1,\frac{y}{x}\right) &= \frac{\sin(\frac{\beta \pi}{2})}{\pi} \frac{1}{x \left(1 - \frac{y}{x} \right)} \left(\frac{1}{- \frac{y}{x}}\right)^\frac{\beta}{2} \One_{(-\infty,0]}\left(\frac{y}{x}\right)\\
 &= \frac{\sin(\frac{\beta \pi}{2})}{\pi} \frac{1}{x-y} \left(-\frac{x}{y}\right)^\frac{\beta}{2} \One_{(-\infty,0]}(y)\\
 & = v_\beta(x,y).
\end{split}
\end{equation*}
\item Using i) yields for $s\in\R$ with substitution $\tilde{y} x = y$
\begin{equation*}
\begin{split}
\int_{-\infty}^s \int_{-\infty}^\infty f(x) v_\beta(x,y)\ dx\ dy & = \int_{-\infty}^\infty \int_{-\infty}^s f(x) \frac{1}{x} v_\beta(1,\frac{y}{x})\ dy\ dx\\
& = \int_{-\infty}^\infty \int_{-\infty}^\infty \One_{(-\infty,s]}(\tilde{y}x) f(x) \frac{1}{x} v_\beta(1,\tilde{y})\,x\ d\tilde{y}\ dx\\
& = \Prob(F V \leq s)
\end{split}
\end{equation*}
and with substitution $-\tilde{y} x = y$
\begin{equation*}
\begin{split}
\int_{-\infty}^s \int_{-\infty}^\infty g(x) u_\alpha(x,y)\ dx\ dy & = \int_{-\infty}^\infty \int_{-\infty}^s g(x) \left(-\frac{1}{x}\right) u_\alpha(-1,-\frac{y}{x})\ dy\ dx\\
& = \int_{-\infty}^\infty \int_{-\infty}^\infty \One_{(-\infty,s]}(-\tilde{y}x) g(x) \left(-\frac{1}{x}\right) u_\alpha(-1,\tilde{y})\,(-x)\ d\tilde{y}\ dx\\
& = \Prob(-G U \leq s).
\end{split}
\end{equation*}
\item Note that
$$\int_0^\infty (y+1)^{-1} y^{-s}\ dy = B(1-s,s) = \frac{\Gamma(1-s)\Gamma(s)}{\Gamma(1)} = \frac{\pi}{\sin\left(s\pi\right)}\ \ \text{ for all } 0< s < 1$$
where $B(\cdot,\cdot)$ is the Beta function and the last equality holds by the reflection formula for the Gamma function (e.g. 6.1.17 in \cite{Abra1972}).
Thus 
$$\E(U^r) = \int_0^\infty y^r u_\alpha(-1,y)\ dy = \frac{\sin\left(\frac{\alpha\pi}{2}\right)}{\pi} \int_0^\infty (y+1)^{-1} y^{-\frac{\alpha}{2} +r}\ dy =\frac{\sin\left(\frac{\alpha \pi}{2}\right)}{\sin\left(\frac{(\alpha-2r)\pi}{2}\right)}$$
for all $r$ such that $\frac{\alpha}{2} - 1 < r < \frac{\alpha}{2}.$ Further for $r\geq \frac{\alpha}{2}$ and $y\geq 1$:
$$(y+1)^{-1} y^{-\frac{\alpha}{2} +r} \geq \frac{1}{2} y^{-\frac{\alpha}{2} +r-1}$$
and this is not integrable on $[1,\infty)$ thus for $r\geq \frac{\alpha}{2}$ the moment is $\infty$. Similarly for $r\leq \frac{\alpha}{2}-1$ and $y\leq 1$:
$$(y+1)^{-1} y^{-\frac{\alpha}{2} +r} \geq y^{-\frac{\alpha}{2} +r}$$
and this is not integrable on $(0,1]$ thus for $r\leq \frac{\alpha}{2}-1$ the moment is $\infty$.

Furthermore for $y>0$
$$v_\beta(1,-y) = \const{\beta} \frac{1}{y+1} y^{-\frac{\beta}{2}} \One_{[0,\infty)}(y) = u_\beta(-1,y)$$
and thus
$$\E((-V)^r) = \int_{-\infty}^\infty (-y)^r v_\beta(1,y)\ dy = \int_{-\infty}^\infty \tilde{y}^r u_\beta(-1,\tilde{y})\ d\tilde{y}$$
and the independence of $V,$ $U$ yields
$$\E((-VU)^r) = \frac{\sin\left(\frac{\alpha \pi}{2}\right)\sin\left(\frac{\beta \pi}{2}\right)}{\sin\left(\frac{(\alpha - 2r) \pi}{2}\right) \sin\left(\frac{(\beta-2r) \pi}{2}\right)}$$
for $r$ in $\left(\frac{\alpha\lor\beta}{2}-1, \frac{\alpha\land\beta}{2}\right)$.

\item For $r^\star = \frac{\alpha+\beta}{4}-\frac{1}{2}$
\begin{equation*}
\begin{split}
\E((-VU)^{r^\star})&=\frac{\sin\left(\frac{\alpha \pi}{2}\right)\sin\left(\frac{\beta \pi}{2}\right)}{\sin\left(\frac{\alpha-\beta}{4}\pi + \frac{\pi}{2}\right) \sin\left(\frac{\beta-\alpha}{4}\pi + \frac{\pi}{2}\right)}\\
&=\frac{\sin\left(\frac{\alpha \pi}{2}\right)\sin\left(\frac{\beta \pi}{2}\right)}{\cos\left(\frac{\alpha-\beta}{4}\pi\right)^2}
=1- \frac{1+\cos\left(\frac{\alpha+\beta}{2}\pi\right)}{1+\cos\left(\frac{\alpha-\beta}{2}\pi\right)}
\end{split}
\end{equation*}
where we used first the translation and symmetry of $\sin$ and $\cos$. In the last step formula 4.3.31 \cite{Abra1972} was used for for the numerator and 4.3.25 \cite{Abra1972} for the denominator.\\
Thus the $r^\star$-moment is less than one for $\alpha+\beta \neq 2.$ Note that $r^\star$ is negative for $\alpha+\beta <2$ and positive for $\alpha+ \beta >2$.
Finally 
\begin{equation*}
\begin{split}
\Prob(-UV \leq s) & = \int \int \One_{(-\infty,s]} (-\tilde{x}\tilde{y})\, v_\beta(1,\tilde{x})u_\alpha(-1,\tilde{y})\ d\tilde{y}\ d\tilde{x}\\
& = \int \int \One_{(-\infty,s]}\left(\frac{-1}{xy}\right)\, v_\beta(1,-\frac{1}{x})u_\alpha(-1,\frac{1}{y})\,\frac{1}{x^2y^2} dy\ dx\\
& = \int \int \One_{(-\infty,s]}\left(\frac{-1}{xy}\right) \const{\alpha}\const{\beta} \frac{1}{1+\frac{1}{x}} x^{\frac{\beta}{2}}\frac{1}{1-\frac{1}{y}} (-y)^{\frac{\alpha}{2}} \frac{1}{x^2y^2} dy\ dx\\
& = \int \int \One_{(-\infty,s]}\left(\frac{-1}{xy}\right) \const{\alpha}\const{\beta} \frac{1}{x+1} x^{\frac{\beta}{2}-1}\frac{1}{y-1} (-y)^{\frac{\alpha}{2}-1} dy\ dx\\
& = \int \int \One_{(-\infty,s]}\left(\frac{-1}{xy}\right) v_\beta(1,y) u_\alpha(-1,x) x^{\frac{\beta+\alpha}{2}-1} (-y)^{\frac{\alpha+\beta}{2}-1} dy\ dx\\
& = \Prob(-(UV)^{-1} \leq s)
\end{split}
\end{equation*}
where we used in the second line the substitution $\tilde{x} = -\frac{1}{x}$ and $\tilde{y} = -\frac{1}{y}$ and for the last step the assumption $\alpha+\beta = 2$.  \qedhere
\end{enumerate}
\end{proof}

\begin{Thm} \label{knotted-stable-thm}
Let $(X_t)_{t\geq 0}$ be a c\`adl\`ag time homogeneous strong Markov process on $\R$ such that \eqref{assume-oscillate} holds and such that there exist $b\in\R,\, \alpha,\beta\in (0,2)$ such that
$$\lim_{t\to 0} \E_x\left(\frac{e^{iX_t\xi}-1}{t} \right) = \begin{cases}
-|\xi|^\beta &\text{ for } x> b,\\
-|\xi|^\alpha &\text{ for } x< b.
\end{cases}$$ 
Then
\begin{enumerate}[\upshape i)]
\item $b$ is left limit recurrent if $\alpha+\beta >2$,
\item $b$ is recurrent if $\alpha+\beta \geq 2$,
\item $b$ is transient if $\alpha+\beta <2$.
\end{enumerate}
\end{Thm}
\begin{proof}
If $b\neq 0$ consider $(X_t-b)_{t\geq 0}$ for which the properties at $0$ correspond to those of $(X_t)_{t\geq 0}$ at $b$. Thus, without loss of generality, we may assume $b=0$.

Let $(Y_n)_{n\geq 0}$ be the overshoot Markov chain corresponding to $(X_t)_{t\geq 0}$ as defined in Section \ref{framework}. Then for $x>0$ and $s\in\R$
$$\Prob_x(Y_n \leq s) = \int_{-\infty}^s \int_{-\infty}^\infty \int_{-\infty}^\infty v_\beta(y,v) u_\alpha(v,u)\ \Prob_x(Y_{n-1}\in dy)\ dv\ du$$
and by Lemma \ref{stable-overshoots-lem} ii) 
$$Y_n \stackrel{d}{=} Y_1 \prod_{i=1}^{n-1} (-U_iV_i)$$
where $U_i \sim u_\alpha(-1,\cdot)$, $V_i\sim v_\beta(1,\cdot)$ and $(U_i)_{i=1,\ldots,n-1},$ $(V_i)_{i=1,\ldots,n-1}$, $Y_1$ are independent. In particular for $r\in \R$
$$\E_x( Y_n^r) = \E_x(Y_1^r) \left(\E\left((-U_1V_1)^r\right)\right)^{n-1}$$
holds and furthermore using the definition of $Y_1$ and Lemma \ref{stable-overshoots-lem} ii) for $\tilde V \sim v_\beta(x,\cdot)$ independent of $U_1$
$$\E_x(Y_1^r) = \E(-\tilde{V}^r) \E(U_1^r) $$
and
$$\E(-\tilde{V}^r) = - \int_\R \tilde{v}^r\  v_\beta(x,\tilde{v})\ d\tilde{v} = - \int_{-\infty}^\infty v^r \frac{1}{x} v_\beta\left(1,\frac{v}{x}\right)\ dv = - x^r \int_\R v^r v_\beta(1,v)\ dv = x^r \E(-V_1^r).$$

To prove i), let $\alpha+\beta >2$ and choose $r>0$, cf. Lemma \ref{stable-overshoots-lem} iv), such that
$$\E((-U_1V_1)^r)< 1.$$
Then $\E_x(Y_1^r) < \infty$ and for all $\varepsilon>0$ by the Chebychev inequality
$$\sum_{n=1}^\infty \Prob_x(Y_n \geq \varepsilon) \leq \sum_{n=1}^\infty \frac{\E_x(Y_n^r)}{\varepsilon^r} = \frac{\E_x (Y_1^r)}{\varepsilon^r} \sum_{n=1}^\infty \E((-U_1 V_1)^r)^{n-1} = \frac{\E_x Y_1^r}{\varepsilon^r} \frac{1}{1-\E((-U_1 V_1)^r)}<\infty$$
holds. Thus the Borel-Cantelli Lemma implies that $Y_n \xrightarrow{n\to\infty} 0$ almost surely. 
Let $q\in (\frac{\alpha \lor \beta}{2}-1,0),$ then $0<\E((-U_1V_1)^q)<\infty$ by Lemma \ref{stable-overshoots-lem} iii). With $R':=\left(2\E((-U_1V_1)^q)\right)^{\frac{1}{q}}$ we get
$$\sup_{y\geq 1} \Prob_y(X_{\sigma_1} < R') = \sup_{y\geq 1} \Prob_y(Y_1 < R') = \sup_{y\geq 1} \Prob_y(Y_1^q > R'^q)\leq \sup_{y\geq 1}\frac{1}{y^{|q|}} \frac{\E((-U_1V_1)^q)}{R'^q} = \frac{1}{2},$$
i.e. \eqref{leftlimit-recurrence-estimate} holds. Thus 0 is left limit recurrent by Theorem \ref{recurrence-transience-thm} iii).

Analogously to prove iii), let $\alpha+\beta <2$ and choose $r<0$, cf. Lemma \ref{stable-overshoots-lem} iv), such that
$$\E((-U_1V_1)^r)< 1.$$
Then $\E_x(Y_1^r) < \infty$ and for all $\varepsilon>0$ by the Chebychev inequality
\begin{equation*}
\begin{split}
\sum_{n=1}^\infty \Prob_x\left(\frac{1}{Y_n} \geq \varepsilon\right) &\leq \sum_{n=1}^\infty \frac{\E_x(Y_n^{-|r|})}{\varepsilon^{|r|}} \\
&= \frac{\E_x (Y_1^r)}{\varepsilon^{|r|}} \sum_{n=1}^\infty \E((-U_1 V_1)^r)^{n-1} = \frac{\E_x Y_1^r}{\varepsilon^{|r|}} \frac{1}{1-\E((-U_1 V_1)^r)}<\infty
\end{split}
\end{equation*}
 holds. Thus the Borel-Cantelli Lemma implies $^{1}/_{Y_n} \xrightarrow{n\to\infty} 0$ almost surely, i.e. $Y_n \xrightarrow{n\to\infty} \infty$ almost surely. 
Now let $q\in \left(0,\frac{\alpha\land\beta}{2}\right),$ then $0<\E((-U_1V_1)^q)<\infty$ and $R:=\left(2\E((-U_1V_1)^q)\right)^{\frac{1}{q}}$ yields 
$$\sup_{y\in(0,1]} \Prob_y(X_{\sigma_1} > R) =\sup_{y\in(0,1]} \Prob_y(Y_1 > R) \leq \sup_{y\in(0,1]} \frac{ \E_y(Y_1^q)}{R^{q}} = \sup_{y\in(0,1]} y^{q}\frac{ \E((-U_1V_1)^q)}{R^{q}} = \frac{1}{2}.$$
Moreover, for $y<0$
$$\Prob_y(X_{\sigma_1} > R) = \int_R^\infty u_\alpha(y,z)\ dz = \int_R^\infty \frac{-1}{y} u_\alpha(-1, -\frac{z}{y}) \ dz = \int^\infty_{-\frac{R}{y}} u_\alpha(-1,\tilde{z})\ d\tilde{z}$$
holds and thus
$$\sup_{y\in[-1,0)} \Prob_y(X_{\sigma_1} > R) = \int_R^\infty u_\alpha(-1,\tilde{z}) \ d\tilde{z} < 1$$
which is strictly less than 1 since $R>0$ and $u_\alpha$ is a probability density with $u_\alpha(-1,x) > 0$ for all $x> 0$. 
Therefore \eqref{local-transience-estimate} holds and 0 is by Theorem \ref{recurrence-transience-thm} i) locally transient. 

Finally let $\alpha+\beta = 2$ and note that
$$\log Y_n \stackrel{d}{=} \log Y_1 + \sum_{i=1}^{n-1} \log(-U_i V_i)$$
holds. By Lemma \ref{stable-overshoots-lem} iv) for any $r\in \R$
$$\Prob(\log(-U_1V_1) \leq r) = \Prob(\log((-U_1V_1)^{-1})\leq r) = \Prob(\log(-U_1V_1) \geq -r)$$
and thus $\log Y_n$ has the same distribution as a symmetric random walk with initial distribution given by $\log Y_1$. Hence
$$\limsup_{n\to\infty} \log(Y_n) = +\infty \text{ and } \liminf_{n\to\infty} \log(Y_n) = -\infty$$
holds and therefore
$$\limsup_{n\to\infty} Y_n = +\infty \text{ and } \liminf_{n\to\infty} Y_n = 0.$$
Now Theorem \ref{recurrence-transience-thm} ii) implies that 0 is locally recurrent.
\end{proof}

\begin{Rem}
Note that we assumed the existence of the process in Theorem \ref{knotted-stable-thm}. The proof of the existence of such a process (and that it is a $\lambda$-irreducible T-model) is part of ongoing research and will be postponed to a forthcoming paper. This seems reasonable to us, since the existence of the process is related to the question of solving SDEs with discontinuous coefficients and the solution theory for such equations requires tools which go beyond the scope of the present paper.
\end{Rem}

The the next result for symmetric $\alpha$-stable L\'evy processes is well known (e.g. \cite{Sato99}). We just present it with a new proof.  

\begin{Cor} 
Let $(X_t)_{t\geq 0}$ be a symmetric $\alpha$-stable L\'evy process with stability index $\alpha\in (0,2)$, then $(X_t)_{t\geq 0}$ is
\begin{enumerate}[\upshape i)]
\item point recurrent if $\alpha >1$,
\item Harris recurrent if $\alpha \geq 1$,
\item transient if $\alpha < 1$.
\end{enumerate} 
\end{Cor}
\begin{proof}
Just apply Theorem \ref{knotted-stable-thm} for $\alpha = \beta$ and note that $b$ can be chosen arbitrary. Further note that the process is clearly a $\lambda$-irreducible T-model, since it is a $C_b$-Feller process with positive transition density. Thus Theorem \ref{dichotomy-tweedie} yields the recurrence-transience dichotomy.\\
Furthermore Lemma \ref{qlc} is applicable since the process is a Hunt process, i.e. in particular it is quasi-left continuous (e.g. Thm. I.9.4 in \cite{blumenthalget}). 
\end{proof}

The results of Section \ref{rec-trans-I} show that two $\lambda$-irreducible $C_b$-Feller processes have the same recurrence (transience) behavior if they have the same generator outside an arbitrary ball. In particular we get the following Corollary for stable-like processes.   

\begin{Cor} Assume the process in Theorem \ref{knotted-stable-thm} exists and is a $\lambda$-irreducible T-model. Let $(X_t)_{t\geq 0}$ be a stable-like process on $\R$ with symbol $|\xi|^{\alpha(x)}$ and suppose there exists $\alpha, \beta \in (0,2)$ such that for some arbitrary $R>0$
$$\alpha(x) = \alpha \text{ for } x < - R,$$
$$\alpha(x) = \beta \text{ for } x> R,$$
then $(X_t)_{t\geq 0}$ is 
\begin{itemize}
\item Harris recurrent if and only if $\alpha+\beta \geq 2$,
\item transient if and only if $\alpha+\beta <2$.
\end{itemize}
\end{Cor}
\begin{proof}
$X_t$ is $\lambda$-irreducible since it has a transition density with respect to the Lebesgue measure (cf. \cite{Nego1994}) and a T-model, since it is a $C_b$-Feller process by Prop. 6.2 in \cite{Bass1988a}.

The process coincides on $\R\backslash\overline{B_R(0)}$ with the process of Theorem \ref{knotted-stable-thm} and therefore by Corollary \ref{tail-comparison} both processes have the same recurrence/transience behavior. Thus Theorem \ref{knotted-stable-thm} implies the result.
\end{proof}

\textbf{Acknowledgement:} The paper was initiated during a visit to Zagreb funded by DAAD. The author is grateful for discussions on the topic with Ren\'e Schilling and Zoran Vondracek.

%\bibliographystyle{model1b-num-names}
%\bibliography{../../bib/MyBib}

\end{document}